\documentclass[12pt,reqno]{amsart}
\makeatletter
\@addtoreset{equation}{section}         

\newcommand {\C}{\mathbb{C}} 
\newtheorem{thm}{Theorem}[section]

\newtheorem{lem}[thm]{Lemma}
\newtheorem{cor}[thm]{Corollary}

\theoremstyle{definition}

\newtheorem{dfn}[thm]{Definition}
\newtheorem{rem}[thm]{Remark}

\newcommand{\myqed}{\ \hfill \rule{2mm}{2mm}}

\newcommand{\beqr}{\begin{eqnarray*}}
\newcommand{\eeqr}{\end{eqnarray*}}
\newcommand{\bc}{\begin{center}}
\newcommand{\ec}{\end{center}}
\newcommand{\LL}{\mathcal{L}}
\newcommand{\codim}{\mathrm{codim \; }}
\newcommand{\dist}{\mathrm{dist}}

\newcommand{\gm}{\gamma}

\newcommand{\dt}{\delta}

\newcommand{\ld}{\lambda}
\newcommand{\sm}{\sigma}

\newcommand{\Om}{\Omega}

\begin{document}
\title{The stability radius
of Fredholm linear pencils}
\author{C. Badea}
\address{Math\'ematiques, UMR 8524 au CNRS, Universit\'e de Lille I,
  F--59655 Villeneuve d'Ascq, France}
\email{badea@agat.univ-lille1.fr}
\author{M. Mbekhta}
\address{Math\'ematiques, UMR 8524 au CNRS, Universit\'e de Lille I,
  F--59655 Villeneuve d'Ascq, France}
\email{mbekhta@agat.univ-lille1.fr}
\date{}

\maketitle

\begin{abstract}
Let $T$ and $S$ be two bounded linear operators from Banach
spaces
$X$ into $Y$ and suppose that $T$ is Fredholm and $\dim N(T-\lambda S)$ is
constant in a neighborhood of $\lambda = 0$. 
Let $d(T;S)$ be 
the supremum of all $r > 0$ such that $\dim
N(T-\ld S)$ and $\codim R(T-\ld S)$ are constant for all 
$\ld$ with $|\ld | < r$. It is a consequence of more general results due to 
H. Bart and D.C. Lay (1980) that $d(T;S) = \lim_{n\to
\infty}\gm_{n}(T;S)^{1/n}$, where $\gm_{n}(T;S)$
are some non-negative
(extended) real numbers. For $X=Y$ and $S = I$, the identity operator,
we have $\gm_{n}(T;S) = \gm (T^n)$, where $\gm$ is the reduced minimum
modulus. A different representation of the stability radius $d(T;S)$ is
obtained here in terms of the spectral radii of
generalized inverses of $T$. The existence of
generalized resolvents for Fredholm linear pencils is also considered.
\end{abstract}


\section{Introduction}
Let $T$ be a bounded linear operator on a Banach space $X$. Denote 
by $\gamma(T)$ the reduced minimum modulus of $T$. The limit
\begin{equation}\label{eq:2}
\lim_{n\to \infty} \gm (T^n)^{1/n}
\end{equation}
was studied for Fredholm operators $T$ by K.H. F\"{o}rster and M.A.
Kaashoek \cite{FoKa}. If $T$ is Fredholm, they proved that the limit 
(\ref{eq:2}) exists and is equal to the supremum of all $r > 0$ such that 
$\dim N(T-\lambda)$ and 
$\codim R(T-\lambda)$ are constant on $\{\lambda : 0< |\lambda| < r\}$. If $T$ is Fredholm
and $0\in reg(T)$, the generalized resolvent set of $T$, then the limit 
(\ref{eq:2}) coincides \cite{Mbe2} with $\dist (0,\sigma_g(T))$, 
the distance of $0$ to the generalized spectrum of $T$. We refer to \cite{FoKa},
\cite{Mbe2} and section 2 of this
paper for the definition of unknown terms and to 
\cite{Apo,BaLa,KoMu,Poo,Sch,Zem1}
for other papers related to the limit
(\ref{eq:2}).

An extension of the limit (\ref{eq:2}) for operator pencils $\ld \to T-\ld 
S$ between two Banach spaces $X$ and $Y$  
has been considered by H. Bart and D.C. Lay \cite{BaLa}. They defined a
sequence of non-negative (extended) real numbers $\gm_{n}(T;S)$
and studied the limit
\begin{equation}\label{eq:3}
\lim_{n\to \infty} \gm_{n}(T;S)^{1/n}\; .
\end{equation}
It was proved in \cite{BaLa} that 
if $T$ is Fredholm, 
then the limit (\ref{eq:3})
exists and coincides with the so-called stability radius $d(T;S)$ 
of $T$ and $S$. In the case when 
$T$ is Fredholm and $\dim N(T-\lambda S)$ is
constant in a neighborhood of $\lambda = 0$, the stability radius $d(T;S)$ 
equals the supremum of all $r > 0$ such that $\dim N(T-\ld
S)$ and $\codim R(T-\ld
S)$ are constant for $|\ld | < r$. 

A different representation formula was given recently \cite{BaMb1}
for the limit (\ref{eq:2}). It was proved there that if $T$ is Fredholm
and $0 \in reg(T)$, the generalized resolvent set of $T$, then
\begin{equation}\label{eq:4}
\lim_{n\to \infty} \gm (T^n)^{1/n} = \sup \{
\frac{1}{r(L)} : TLT = T\}\; ,
\end{equation}
where $r(L)$ is the spectral radius of $L$. 
The same result (\ref{eq:4}) 
is true without the condition of Fredholmness of $T$
for Hilbert space operators \cite{BaMb2}. 

The aim of the present note is to extend to linear pencils the formula
(\ref{eq:4}). Namely, we
will prove that
\begin{eqnarray*}
	d(T;S) & = & \sup \{\frac{1}{r(SL)} : L \in \LL (Y,X), \;  
	TLT = T\}  \\
	 & = & \sup \{\frac{1}{r(SL)} : L \in \LL (Y,X), \;  
	 TLT = T ,  LTL = L\} ,
\end{eqnarray*}
if $T$ is Fredholm and $\dim N(T-\lambda S)$ is
constant in a neighborhood of $\lambda = 0$. The above formula
relates the limit (\ref{eq:3}) with the spectral radii of
generalized inverses of $T$. It is an open problem if this representation 
holds without the condition of
constancy of $\dim N(T-\lambda S)$.

One of the ingredients of the proof of the main result is the
existence
of generalized resolvents for some operator pencils. 
This is
interesting in its own. 

The paper is organized as follows. In the next section we recall some 
notation and known results. In section 3 we introduce and study 
generalized resolvents for operator pencils and their existence 
for some operator pencils is proved in section 
4. This is used in section 5 to prove the main result concerning 
the stability radius.
\section{Notation and known results}
We present in this section some notation, basic definitions and known 
results. 

Recall that $X$ and $Y$ will denote complex Banach spaces.
We denote by $\LL(X,Y)$ the
Banach space of all continuous, linear operators from $X$ into $Y$~; 
we abbreviate $\LL(X,X)$ to $\LL(X)$. We
use $R(T) = TX$, $N(T) = \{x \in X : Tx = 0\}$, $\sm (T)$ and
$\rho (T) = \C \setminus \sm (T)$ to denote
the range, the
kernel, the spectrum and, respectively, the resolvent set of $T$. We
denote by $r(T)$ the spectral radius of $T$. 
We write $X = E \oplus F$ to designate
direct sums (i.e. $E \cap F = \{ 0\}$ and $E+F = X$). 
In this case we say that $F$ is the (direct)
complement of $E$ in $X$. The distance from a point $x$ to a set $A$ 
is denoted by $\dist (x,A)$ and $B(\ld_{0},r)$ is the open set  
$\{\ld \in \C : |\ld - \ld_{0}| < r\}$.

The operator $T \in \LL (X,Y)$ is \emph{Fredholm} if 
both the dimension of $N(T)$ and the codimension of
$R(T)$ are finite. The range of a Fredholm operator
is closed.

The
\emph{reduced minimum modulus} of $T \neq 0$, $T \in \LL (X,Y)$,
is defined by 
$$\gm (T) = \inf \{\| Tx\| : x \in X , \dist (x, N(T)) = 1\}.$$
If $T = 0$ we set $\gm (T) = \infty$.

Let $T \in \LL(X,Y)$. An operator $L \in \LL(Y,X)$ is called a 
\emph{generalized inverse} of $T$ if $TLT = T$ and $LTL = L$. If $L$ 
satisfies only the first condition $TLT = T$, then $L_{1} = LTL$ is a 
generalized inverse of $T$.

Let $T \in \LL (X)$. Let $reg(T)$ denote \cite{Mbe2},  \cite{Mbe3} 
the \emph{generalized resolvent 
set} of $T$, that is the set of all complex numbers
$\ld$ for which
there exists a neighborhood $V_{\ld}$ of $\ld$ and an analytic 
function $G$ on $V_{\ld}$ such that $G(\mu)$ is a generalized 
inverse of $T - \mu I_{X}$ for each $\mu \in V_{\ld}$. 
Then (cf. \cite{Mbe2,Mbe3}) $0\in reg(T)$ if and only if $T$ has a 
generalized inverse
and $N(T) \subseteq R(T^m)$ for every $m$. Let 
$\sigma_g(T) = \C \setminus reg(T)$ be the \emph{generalized spectrum} of $T$.
 
Let $T,S \in \LL (X,Y)$ and consider the linear operator pencil
$\ld \to T - \ld S \in \LL (X,Y)$.
Define the subspaces $N_{m} = N_{m}(T;S)$ and $R_{m} = R_{m}(T;S)$ of
$X$ by $N_{0} = \{ 0\}, R_{0} = X$ and 
$$N_{m+1} = T^{-1}SN_{m} \quad , \quad R_{m+1} = S^{-1}TR_{m}.$$
for $m = 0, 1, \cdots$. Note that the sequence $(R_{m})$ is
decreasing. Define
$$X_{\infty} = \bigcap_{m} R_{m} \quad ; \quad Y_{\infty} = \bigcap_{m}
TR_{m} \; .$$
The extended integer $k(T;S)$ defined by 
$$k(T;S) = \dim N(T)/\left( N(T)\cap X_{\infty}\right)$$
is the \emph{stability number} of $T$ and $S$. So 
$k(T;S) = 0$ means $N(T) \subset R_{m}$ for all $m$. 
For further reference we note that $k(T;S) = 0$ implies
(cf. \cite[Lemma 2.3]{Kaa}) $TX_{\infty} = Y_{\infty}$ and
$S^{-1}Y_{\infty} = X_{\infty}$. If $T$ is Fredholm, then $k(T;S) = 0$ if and only if 
$\dim N(T-\lambda S)$ is
constant in a neighborhood of $\lambda = 0$, if and only if $\codim R(T-\lambda S)$ is
constant in a neighborhood of $\lambda = 0$. We will use these equivalences several times.

Let $m \geq 1$. The $m$-tuple $(x_{1},x_{2}, \cdots , x_{m}) \in X^m$ is a
\emph{chain} for $T$ and $S$ if $Tx_{i} = Sx_{i-1}$ for 
$i = 2,\cdots , m$.
Let $\gm_{m} = \gm_{m}(T;S)$ denote the supremum of all $c \geq 0$
with the property that $\| Tx_{1}\| \geq c \; \dist (x_{m},N_{m})$
for all chains $(x_{1},x_{2}, \cdots , x_{m})$. Then $\gm_{1}(T;S) = 
\gm(T)$
is the reduced minimum modulus of $T$.

Let $T,S \in \LL(X,Y)$. Suppose that $R(T)$ is closed and $k(T;S) < 
+\infty$. 
The \emph{stability radius} $d(T;S)$ of $T$ and $S$ 
is defined \cite{BaLa} 
as the supremum of 
all $r > 0$ such that $R(T-\ld S)$ is closed and $k(T-\ld S) = 0$ for 
$0 < |\ld | < r$. When $T$ is Fredholm, the stability radius $d(T;S)$ is 
equal to the supremum of all $r > 0$ such that $\dim N(T - \ld S)$ 
and $\codim R(T-\ld S)$ are constant on $0 < |\ld | < r$. 
If $T$ is Fredholm, 
then \cite{BaLa} $d(T;S) = \lim_{m\to \infty} \gm_{m}(T;S)^{1/m}$.
\section{Generalized resolvents for linear pencils}
We start with the following definition. 
\begin{dfn}
		Let $T$ and $S$ be two elements of $\LL (X,Y)$. 
		Let $U$ be an open set in the complex plane. The function 
$$U \ni \ld \to G(\ld) \in \LL(Y,X)$$
is called a \emph{generalized resolvent on} $U$ of
the linear pencil $\ld \to T-\ld S \in \LL(X,Y)$ if
\begin{itemize}
			\item[(1)]  $(T - \ld S)G(\ld)(T - \ld S) = T - \ld S$,
for all $\ld \in U$ ;

			\item[(2)]  $G(\ld)(T - \ld S)G(\ld) = G(\ld)$, for all
$\ld \in U$ ;

			\item [(3)]  $G(\ld) - G(\mu) = (\ld - \mu)G(\ld)SG(\mu)$
for all $\ld$ and $\mu$ in the \emph{same} connected component of $U$.
		\end{itemize}
\end{dfn}

The first two conditions say that $G(\ld)$ is a generalized inverse 
of $T-\ld S$ for each $\ld$, while the third one is an analogue of the 
classical resolvent identity. 
The assumption that 
the resolvent identity (3) holds only for $\ld$ and $\mu$ in 
the same connected component is consistent with \cite{Mbe3}. It is 
possible to have generalized resolvents with  
(3) not fulfilled for $\ld$ and $\mu$ in distinct connected 
components \cite[p. 376]{Mbe3}.
We refer to \cite{AC, Mbe2,LM, Mbe3,BaMb3} 
for properties of 
generalized resolvents in the classical case $X = Y$, $S = I_{X}$. 

The following lemma shows that each generalized resolvent on $U$ of a 
linear pencil is analytic on $U$.

\begin{lem}\label{lem:32}
	Let $U$ be an open set in the complex plane and let 
	$$U \ni \ld \to G(\ld) 
	\in \LL(Y,X)$$ 
	be an operator function satisfying 
	$$G(\ld) - G(\mu) = (\ld - \mu)G(\ld)SG(\mu)$$
for all $\ld$ and $\mu$ in the same connected component of $U$. 
Then $G$ is analytic on $U$.
\end{lem}

\noindent {\bf Proof .} Let $\ld_{0}$ be a fixed point in $U$ and 
let $r$ be a positive number such that $r < \|SG(\ld_{0})\|^{-1}$ and 
$B(\ld_{0},r)$ is included in the connected 
component of $U$ containing $\ld_{0}$. Let $\ld \in B(\ld_{0},r)$. 
The resolvent 
identity (3) 
implies that
$$G(\ld_{0}) = G(\ld)\left[ I - (\ld - \ld_{0})SG(\ld_{0})\right]$$
and thus
\begin{eqnarray*}
	G(\ld) & = & G(\ld_{0})\left[ I - (\ld - \ld_{0})SG(\ld_{0})\right]^{-1}  \\
	 & = & G(\ld_{0})\sum_{k\geq 0}(\ld - \ld_{0})^k (SG(\ld_{0}))^k 
\end{eqnarray*}
for each $\ld \in B(\ld_{0},r)$. 
Thus $G$ is analytic on $U$. \myqed

\begin{dfn}
	Let $U$ be an open set in $\C$. 
	The linear pencil $\ld \to T-\ld S \in \LL(X,Y)$
	is said to have \emph{fixed complements on} $U$ if for each connected 
	component $\Gamma$ of $U$ there exist
	two closed subspaces $E$ and $F$
	of $X$ and $Y$ such that $X = N(T-\ld S) \oplus E$ and $Y = R(T-\ld S)
	\oplus F$
	for all $\ld \in \Gamma$.
\end{dfn}

\begin{rem}\label{rem:thij}
	Consider a linear pencil $\ld \to T-\ld S \in \LL(X,Y)$ and let $w$ 
	be a complex number such that $R(T - wS)$ is closed and there exists 
	a bounded 
	generalized inverse for $T - wS$. 
	It follows from \cite[Proposition I.2.2]{thij} that there 
	exists a (connected) neighborhood $V_{w}$ of $w$ such that 
	$\ld \to T-\ld S$ has fixed complements on $V_{w}$ if and only if the linear 
	pencil $\ld \to T-\ld S$ is uniformly regular at $w$
	\cite{thij}, that is the function $\ld \to \gm(T-\ld S)$ is continuous 
	at $w$.
\end{rem}

The following result gives a (global) characterization of the 
property of having fixed complements in terms of generalized 
resolvents.  

\begin{thm} \label{fixed}
Let $T$ and $S$ be two elements of $\LL (X,Y)$.  
	Let $U \subset \C$ be an open set. 
	There exists a generalized resolvent for $\ld \to T - \ld  S$ 
	on $U$ if and only if 
the linear pencil $\ld \to T - \ld S$ has fixed
complements on $U$. 
\end{thm}

\noindent {\bf Proof.} Via reduction to connected components, it is 
possible to assume that $U$ is connected.

Suppose that $\ld \to T - \ld S$ has fixed
complements
on $U$ ; that is, there exist
two closed subspaces $E$ and $F$
of $X$ and $Y$ such that $X = N(T - \ld S) \oplus E$ and $Y = R(T - \ld S)
\oplus F$
for all $\ld \in U$. We will consider the projections $P(\ld)$ and
$Q(\ld)$ onto $R(T - \ld S)$ along
		$F$ and, respectively, onto $E$ along $N(T-\ld S)$.

We will prove the existence of an operator valued function
$G(\ld) \in \LL (Y,X)$, $\ld \in U$, such that
		\begin{itemize}
			\item[(1)]  $(T - \ld S)G(\ld)(T - \ld S) = T - \ld S$,
for all $\ld \in U$ ;

			\item[(2)]  $G(\ld)(T - \ld S)G(\ld) = G(\ld)$, for all
$\ld \in U$ ;

			\item [(3)]  $G(\ld) - G(\mu) = (\ld - \mu)G(\ld)SG(\mu)$
for all $\ld$ and $\mu$ in $U$.
		\end{itemize}
The generalized resolvent $G$ will be analytic on $U$ by Lemma \ref{lem:32}.

Let $u \in Y$. Then $P(\lambda)u \in
R(T - \ld S)$.  Therefore, there exists $v \in X$
such that $P(\lambda)u = (T - \ld S)v$. Set $G(\lambda)u =
Q(\lambda)v$. 

This definition is correct. Indeed, if $w \in X$ is such
that $(T - \ld S)w = P(\lambda)u = (T - \ld S)v$, then $v - w \in N(T -
\ld S)
= N(Q(\lambda))$. Therefore $Q(\lambda)v = Q(\lambda)w$.
Hence $G(\lambda)$ does not depend on the choice of $v$.

Note that $G(\lambda)$ is the generalized inverse of $T- \lambda S$ with null space $F$ and range $E$.
This gives the proofs of (1) and (2).

We also have 
\begin{equation}
	(T - \ld S)G(\ld) = P(\ld) 
	\label{eq3}
\end{equation}
and 
\begin{equation}
	G(\ld)(T - \ld S) = Q(\ld) \; .
	\label{eq2}
\end{equation}

Using these two equations and the fact that 
$G(\ld)$ is a generalized inverse of $T - \ld S$, we get
\begin{equation}
	Q(\ld)G(\ld) = G(\ld)
	\label{b1}
\end{equation}
and
\begin{equation}
	G(\ld)P(\ld) = G(\ld) \; .
	\label{b2}
\end{equation}
The projection $I - P(\ld)$ is onto $F$ and
thus $P(\mu)(I - P(\ld)) = 0$. Similarly, $(I - Q(\ld))Q(\mu) = 0$.
We obtain
\begin{equation}
	P(\mu)P(\ld) = P(\mu) 
	\label{a1}
\end{equation}
and

\begin{equation}
	Q(\ld)Q(\mu) = Q(\mu).
	\label{a2}
\end{equation}
Now we can write
\begin{eqnarray*}
	G(\ld) - G(\mu) & = & G(\ld)P(\ld)-Q(\mu)G(\mu)  \textrm{ (using
(\ref{b1}) and (\ref{b2})) }\\
	 & = & G(\ld)P(\ld)P(\mu)-Q(\ld)Q(\mu)G(\mu) \\
	 & & \textrm{ (using (\ref{a1})
	 and (\ref{a2})) }  \\
	 & = & G(\ld)P(\mu)-Q(\ld)G(\mu) \textrm{ (using (\ref{b1}) and
	 (\ref{b2})) }  \\
	 & = & G(\ld)(T - \mu S)G(\mu)-G(\ld)(T - \ld S)G(\mu) \\
	 & & \textrm{ (using
	 (\ref{eq2}) and
	 (\ref{eq3})) }  \\
	 & = & (\ld - \mu)G(\ld)SG(\mu).
\end{eqnarray*}
This proves that $G(\ld)$ is a generalized resolvent of the linear
pencil $\ld \to T - \ld S$.

For the second part, suppose a generalized resolvent $G\in \LL(Y,X)$ exists
on $U$. We will show the existence of the fixed complements $E$ and
$F$ of $\ld \to T - \ld S$.

\noindent {\scshape STEP} a \quad \emph{The subspace $E$}

Using the resolvent identity, we have
$
G(\ld) - G(\mu) = - [G(\mu) - G(\ld)] = (\ld - \mu)G(\mu)SG(\ld)
$.
We infer that $G(\ld) = G(\mu) + G(\mu)(\ld - \mu)SG(\ld)$,
yielding $R(G(\ld)) \subseteq R(G(\mu))$. Reversing the role of
$\ld$ and $\mu$ we get $R(G(\mu)) = R(G(\ld))$. Take $E = R(G(\ld)) 
\subset X$. Note that the range of $G(\lambda)$, the generalized inverse 
of $T-\lambda S$, is closed.

\noindent {\scshape STEP} b \quad \emph{The subspace $F$}

The resolvent identity (3) implies $N(G(\mu) \subseteq N(G(\ld))$ and 
thus $N(G(\mu) = N(G(\ld))$. We denote $F = N(G(\ld))$.

\noindent {\scshape STEP} c \quad \emph{$Y = R(T-\ld S) \oplus 
F$}

We have $u = [I - (T-\ld S)G(\ld)]u + (T-\ld S)G(\ld)u$ for each $u \in
Y$.
The first term of the sum is in $N(G(\ld)) = F$ since $G(\ld)$ is a
generalized inverse of $T-\ld S$. The second term of the sum is clearly
in $R(T-\ld S)$.

If $y = (T-\ld S)x$ and $y \in F = N(G(\ld)$, then
$0 = G(\ld)y = G(\ld)(T-\ld S)x$,
yielding $0 = (T-\ld S)G(\ld)(T-\ld S)x = (T-\ld S)x = y$.

\noindent {\scshape STEP} d \quad \emph{$X = N(T-\ld S) \oplus E$}

The proof of this equality is similar to that of Step c and
will be omitted.  \myqed

\section{Linear pencils with fixed complements}
\begin{thm}\label{thm:modified}
	Assume that for each $\ld$ in an open, connected set $\Om \subset \C$ the 
	operator $T-\ld S \in \LL(X,Y)$ is Fredholm, and that the dimension 
	of the null space $N(T-\ld S)$ does not depend on $\ld \in \Om$. Let 
	$U$ be a bounded, open and connected set with $U \subset \overline{U} 
	\subset \Om$. Then there exists a generalized resolvent for $\ld \to 
	T-\ld S$ on $U$.
\end{thm}

Using Theorem \ref{fixed}, we have to show that 
$\ld \to T-\ld S$ has fixed complements on $U$, that is, 
there exist two
closed subspaces $E$ of $X$ and $F$ of $Y$ 
such that $X = N(T-\ld S) \oplus E$ and $Y = R(T-\ld S) \oplus F$
for every $\ld \in U$. The existence of generalized resolvents (and thus 
of fixed complements) has been studied in 
\cite{Sap,AC,Mbe2,LM,Mbe3,BaMb3} in the case $X=Y$, $S=I$. 
Proofs of Theorem \ref{thm:modified} in the case $X=Y$, $S=I$ 
can be found in \cite{Sap} and \cite{Mbe3}.
	
For the proof we will need the following key result. Variants of this result
can be found in \cite{Sap,Zem2,Mbe3}. 
\begin{thm}\label{key}
	Assume that for each $\ld$ in an open, connected set $\Om \subset \C$ the 
	operator $T-\ld S \in \LL(X,Y)$ is Fredholm, and that the dimension 
	of the null space $N(T-\ld S)$ does not depend on $\ld \in \Om$.
	Let $W$ be a closed subspace of $Y$ such that
\begin{equation}\label{eq:star}
\codim [R(T - \ld S)\oplus W] \geq 1, \quad \mbox{ for
all }\ld \in \Om .
\end{equation} 
If $U$ is a bounded, open, connected set such that $U \subset
	\overline{U} \subset \Om$, then there exists $x \in Y$ such that $x
	\not\in R(T - \ld S)\oplus W$ for all $\ld \in U$.
\end{thm}

We start with the following Lemma.
\begin{lem}\label{lem:analytic}
	Assume that for each $\ld$ in an open, connected set $\Om \subset \C$ the 
	operator $T-\ld S \in \LL(X,Y)$ is Fredholm, and that the dimension 
	of the null space $N(T-\ld S)$ does not depend on $\ld \in \Om$. 
	Let $W$ be a closed subspace of $Y$. Then there exists an analytic function 
	$$\Om \ni \ld \to P(\ld) \in \LL(Y\oplus W)$$ 
	such that 
	\begin{equation}
		N(P(\ld)) = R(T - \ld S)\oplus W , \quad \ld \in \Om.
		\label{eq:41}
	\end{equation} 
\end{lem}

\noindent {\bf Proof.} Consider the 
Fredholm operator function 
$$F(\ld) = (T-\ld S)\oplus J : X\oplus W \to Y,$$
where $J : W \to Y$ is the inclusion operator. Then $R(F(\ld))$ has 
constant finite codimension on $\Om$ and $N(F(\ld)) = N(T-\ld S)$ has 
constant finite dimension. Then (cf. for instance \cite{BKL}), 
the function $F(\ld)$ has a 
global analytic generalized inverse $G(\ld)\in \LL(Y,X\oplus W)$ on 
$\Om$. 
Take 
$$P(\ld) = I_{X\oplus W} - F(\ld)G(\ld) \in \LL(X\oplus W).$$ 
Then $N(P(\ld)) = R(F(\ld)) = R(T-\ld S) \oplus W$. \myqed

\noindent {\bf Proof of Theorem \ref{key}.}
Let $W$ be a closed subspace of $Y$ as in Theorem \ref{key}. 
Using the above Lemma, there exists an analytic function $P(\ld)$ 
such that (\ref{eq:41}) holds. 

Consider $\ld_{0} \in U$ and let $u \in Y$ be
such that $u \notin R(T - \ld_{0}S)\oplus W$.
The existence of such an element $u$ follows from condition 
(\ref{eq:star}). Then $P(\ld_{0})u \neq 0$. 

Two situations can occur. 
If $\ld \to P(\ld)u$ does not vanish in $\overline{U}$, 
then we take $x = u$. If not, the analytic function 
$\ld \to P(\ld)u$ vanishes in the
compact $\overline{U}$ for finitely many points $\ld_{1}, \ld_{2}, \cdots,
\ld_{n}$, each with finite order of multiplicity $d_{1}, \cdots
d_{n}$.

\noindent {\bf Claim.} There exists $y \in Y$ such that
the analytic function $\ld \to P(\ld)y$ vanishes in $\overline{U}$
only for $\ld_{2}, \cdots , \ld_{n}$.

The proof of this claim is deffered
to the next paragraph. Assuming this claim,
we note that continuing this construction we
are sure that in at most $n$ steps we obtain an element $z \in Y$ such
that
$P(\ld)z \neq 0$. Then this $z$ does not belong to $R(T - \ld S)\oplus W$
for any $\ld \in U$.

\noindent {\bf Proof of the Claim.} The analytic function
$\ld \to P(\ld)u$ vanishes in the
compact $\overline{U}$ for $\ld_{1}, \ld_{2}, \cdots,
\ld_{n}$. Since $P(\ld_{1})u = 0$, we infer that $u \in R(T -
\ld_{1}S)\oplus
W$. Thus there exist $v_{1} \in X$ and $w_{1} \in W$ such that
$u = (T - \ld_{1}S)v_{1} + w_{1}$.
For any $\ld \in \overline{U}$, $\ld \neq \ld_{1}$, we have
\begin{eqnarray*}
	P(\ld)u & = & P(\ld)(T - \ld_{1}S)v_{1} + P(\ld)w_{1}  \\
	 & = & P(\ld)(T - \ld S + \ld S - \ld_{1}S)v_{1} 
= (\ld - \ld_{1})P(\ld)Sv_{1}.
\end{eqnarray*}
The equality
$P(\ld)u = (\ld - \ld_{1})P(\ld)Sv_{1}$
 also holds for $\ld = \ld_{1}$ and thus it holds for all $\ld \in
 \overline{U}$. In particular, zeros of $P(\ld)Sv_{1}$ in
 $\overline{U}$ with $\ld \neq \ld_{1}$ are zeros of $P(\ld)u$.

If $P(\ld_{1})Sv_{1} \neq 0$, then the representation of $P(\ld)u$ shows
 that $P(\ld_{1})Sv_{1}$ vanishes in $\overline{U}$ only for
 $\ld_{2}, \cdots , \ld_{n}$ and we take $y = Sv_{1}$.

If $P(\ld_{1})Sv_{1} = 0$, then $Sv_{1} \in R(T - \ld_{1}S)\oplus W$. We
can write $Sv_{1} = (T-\ld_{1}S)v_{2} + w_{2}$, $v_{2} \in X$. As above,
we
have $P(\ld)Sv_{1} = (\ld - \ld_{1})P(\ld)Sv_{2}$. Then
$P(\ld)u = (\ld - \ld_{1})P(\ld)Sv_{1} = (\ld -
\ld_{1})^2P(\ld)Sv_{2}$
for all $\ld \in \overline{U}$. 

If $P(\ld_{1})Sv_{2} \neq 0$, then we can 
take $y = Sv_{2}$. If $P(\ld_{1})Sv_{2} =
0$, then there exists $v_{3}
\in X$ such that
$P(\ld)u = (\ld -
\ld_{1})^3P(\ld)Sv_{3}$.

This construction will lead eventually to an element $y = Sv_{k}$
satisfying the conditions of the Claim
since the order $d_{1}$ of multiplicity of $\ld_{1}$ is finite. This
completes the proof of the claim and of the theorem. \myqed

\bigskip

\noindent {\bf Proof of Theorem \ref{thm:modified}.} 
By considering the adjoint pencil $T^* - \ld S^*$ acting between the 
dual spaces of $Y$ and $X$, it is sufficient to prove the existence of a fixed
complement $F$ of $R(T - \ld S)$ for $\ld \in U$. 

Without loss of any generality we can
assume $0 \in U$. Let $n = \codim R(T)$. If $n = 0$, take $F = \{ 0\}$.
If $n \geq 1$, by Theorem \ref{key} (with $W = \{ 0 \}$),
there exists $x \in Y$ such
that $x \not\in R(T-\ld S)$ for all $\ld \in U$. Then
$R(T-\ld S) \cap \mbox{ lin} (x) = \{ 0\}$ and
$\codim ( R(T-\ld S) \oplus \mbox{ lin} (x) ) = n-1$. Here $\mbox{
lin}(x)$ is the subspace spanned by $x$. Note that $R(T-\ld S) \oplus
\mbox{lin} (x)$ is a closed subspace of $Y$. By a repeated application of
Theorem \ref{key} (with suitable $W$),
we get a closed subspace $F$ such that $R(T-\ld S) \cap
F = \{ 0\}$ and $\codim R(T-\ld S) \oplus
F =  0$. Therefore $ R(T-\ld S) \oplus F = Y$ for all $\ld \in U$. \myqed 

\begin{cor}\label{cor:main}
Assume that $T\in \LL(X,Y)$ is Fredholm and $\dim N(T-\lambda S)$ is
constant in a neighborhood of $\lambda = 0$. Let $d = 
d(T;S)$ be the stability radius. Let 
$U \subset \overline{U} \subset B(0,d)$ be a bounded, open, 
connected subset of the open ball of radius $d$. Then there 
exists a generalized resolvent for $T-\ld S$ on $U$.
\end{cor}

\noindent {\bf Proof.} Recall that $k(T;S) = 0$. By \cite[Lemma 3.2]{BaLa}, 
$T-\ld S$ is Fredholm for each $\ld \in \Om:= B(0,d)$ 
and \cite[p. 309]{BaLa} the dimension 
of the null space $N(T-\ld S)$ does not depend on $\ld \in \Om$. The 
conclusion follows now from Theorems \ref{thm:modified} and \ref{fixed}. 
\myqed
\section{{\bf The stability radius}}
The following result is the announced formula for the stability radius.

\begin{thm}\label{thm:main}
Let $T$ and $S$ be two elements of $\LL (X,Y)$ such that $T$ is Fredholm
and $\dim N(T-\lambda S)$ is
constant in a neighborhood of $\lambda = 0$. 
Let $d(T;S)$ be the stability radius of $T$ and 
$S$, which equals the supremum of all $r > 0$ such that $\dim
N(T-\ld S)$ and $\codim R(T-\ld S)$ are constant for all 
$\ld$ with $|\ld | < r$. Then we have 
\begin{eqnarray*}
	d(T;S) & = & \sup \{\frac{1}{r(SL)} : L \in \LL (Y,X), \;  
	TLT = T\}  \\
	 & = & \sup \{\frac{1}{r(SL)} : L \in \LL (Y,X), \;  
	 TLT = T,  LTL = L\} .
\end{eqnarray*}
\end{thm}

\medskip

\noindent We start with the following auxiliary results. 

\begin{lem}\label{lem:conorme}
Suppose there exists $\delta > 0$ such that
$$(T- \ld S) F(\ld)(T- \ld S) = T- \ld S$$
for every $\ld$ with $|\ld| < \dt$ and $F(\ld) = \sum_{n\geq 0} 
\ld^nF_{n}$ for $|\ld| < \dt$. Then
$$\gm_{m}(T;S) \geq \frac{1}{\|F_{m-1}\|}$$
for every $m\geq 1$.
\end{lem}

\noindent {\bf Proof.} Let $(x_{1}, \ldots , x_{m})$ be a chain for 
$T$ and $S$. Define 
$$\psi (\ld) = \left[ I_{X} - F(\ld)(T-\ld S)\right]\left( 
\sum_{i=1}^m \ld^{i-1}x_{i}\right) .$$
We have $(T-\ld S)\psi(\ld) = 0$. Denote by $\psi_{k}$, $k \geq 0$, 
the Taylor coefficients of $\psi$ around zero. We then have $\psi_{0} 
\in N(T) = N_{1}$ and $T\psi_{k} = S\psi_{k-1}$. We get recursively 
$\psi_{k-1} \in N_{k}$ for each $k \geq 0$. 
We have  
$$\psi_{m-1} = x_{m} - F_{m-1}Tx_{1}$$ 
and therefore $x_{m} - F_{m-1}Tx_{1} \in N_{m}$. This implies 
$$\dist (x_{m},N_{m}) \leq \|F_{m-1}Tx_{1}\| \leq \|F_{m-1}\| \; 
\|Tx_{1}\|.$$
Since this holds for every chain $(x_{1}, \ldots , x_{m})$, we obtain 
the desired estimate for $\gm_{m}(T;S)$.
\myqed

\begin{lem}\label{lem:2}
	Let $T,S\in \LL(X,Y)$ and suppose that $N(T) \subset R_{m}$ for all $m$. 
	Let $L\in \LL(Y,X)$ with $TLT = T$. 
	Set $\alpha = \min (\|SL\|^{-1} ; \|LS\|^{-1} )$ and 
$$F(\ld) = L\left( I - \ld SL\right)^{-1} = \sum_{k\geq 0}L(SL)^k\ld^k$$
for $\ld$ satisfying $|\ld| < \alpha$.
Then 
\begin{equation}
	(T- \ld S) F(\ld)(T- \ld S) = T- \ld S
	\label{eq:last}
\end{equation}
for each $\ld$, $|\ld| < \alpha$.
\end{lem}

\noindent {\bf Proof.} Recall that $N(T) \subset R_{m}$ for all $m$ 
is equivalent to $k(T;S) = 0$. 
We have $I - F(\ld)(T - \ld S) = \sum_{k\geq 0}(LS)^k(I-LT)\ld^k ,$ the last
series being convergent for $|\ld| < \alpha$.
This shows that
$$(T- \ld S) F(\ld)(T- \ld S) = T- \ld S - A(\ld),$$
where $A(\ld) = (T-\ld S)\left[I - F(\ld)(T - \ld S)\right]$
is given by
\begin{equation}\label{eq:52}
A(\ld) = \sum_{k\geq 1}(TL - I)(SL)^{k-1}S(I - LT)\ld^k .
\end{equation}

The proof that all Taylor coefficients of $A(\ld)$ are zero 
is obtained in several steps.

\noindent {\scshape Step} $\alpha$ \quad \emph{We have $S(I - 
LT)x \in Y_{\infty}$ for every $x \in X$.}

Indeed, for every $x \in X$, $(I - LT)x \in N(T)$. 
Since $N(T) \subset R_{m}$ for all $m$, we have $(I - 
LT)x \in X_{\infty}$. Therefore \cite[Lemma 2.3]{Kaa}
$$S(I - LT)x \in SX_{\infty} \subseteq Y_{\infty}.$$

\noindent {\scshape Step} $\beta$ \quad \emph{$L(Y_{\infty}) \subseteq
X_{\infty}$.}

We have \cite[Lemma 2.3]{Kaa} $TX_{\infty} = Y_{\infty}$. 
If $y \in Y_{\infty}$, $y = Tx$, $x \in X_{\infty}$, then
$Ly = LTx = x + (LT - I)x \in X_{\infty}$.
The latter follows from $x \in X_{\infty}$, $(LT - I)x \in N(T)$ and
from $N(T) \subset X_{\infty}$.

\noindent {\scshape Step} $\gamma$ \quad \emph{For any $k \geq 1$, $(SL)^k
(Y_{\infty}) \subseteq Y_{\infty}$.}

Let $y \in Y_{\infty}$. Using
Step $\beta$ we get $Ly \in X_{\infty}$  ; thus $SLy \in
Y_{\infty}$. Applying this $k$ times, we get $(SL)^k
(Y_{\infty}) \subseteq Y_{\infty}$.

\noindent {\scshape Step} $\delta$ \quad \emph{The restriction of $TL$ on
$Y_{\infty}$
acts like the identity operator on $Y_{\infty}$. }

Indeed, if $y \in Y_{\infty}$, then there exists $x \in
X_{\infty}$ such that $y = Tx$. Then
$TLy = TLTx = Tx = y$.

By (\ref{eq:52}) we have $A(\ld) = 0$ and thus 
(\ref {eq:last}) holds
for each $\ld$ satisfying $|\ld| < \alpha$. 
\myqed

\noindent {\bf Proof of Theorem \ref{thm:main}.} 
Note that the first sup is greater or equal than the second.

Let $\Om = B(0,d)$ be the open disk of radius 
$d = d(T;S)$. Let $\varepsilon > 0$.
Consider the disk $U = B(0,d/(1+\varepsilon)$.
Then $U \subset
\overline{U} \subset \Om$. 

It is a consequence of the hypothesis that 
$T$ is Fredholm and $k(T;S) = 0$. 
This gives the case $\lambda = 0$ of the following 
more general statement~:
The operator $T -\ld S$
is Fredholm and $k(T-\ld S; S) = 0$ for each $\ld \in \Om$. For $\lambda \neq 0$, this
follows by combining  
\cite[Lemma 3.2]{BaLa} with \cite[Theorem 3.1]{BaLa}. 

By
Corollary \ref{cor:main}, there exists a generalized resolvent 
$G$ for the pencil $\ld \to T-\ld S$ on $U$.
Set $G(\ld) = \sum_{n=0}^{\infty}\ld^n G_{n}$, $\ld \in U$.
Since $G$ satisfies the resolvent identity on $U$, we have 
$G(\ld) - G_{0} = G(\ld) - G(0) = \ld G(\ld)SG(0)$.
This implies $G_{n} = G_{n-1}SG_{0}$. Therefore
$G_{1} = G_{0}SG_{0}$ and, by recurrence, $G_{n} = G_{0}(SG_{0})^n$
for all $n$. 

Denote $M = \max\{ \| G(\lambda)\| : |\lambda| \leq d/(1+2\varepsilon)\}$. 
Using Cauchy's estimates, we obtain 
$$\| G_{0}(SG_{0})^n \| = \|G_{n}\|
\leq M\left(\frac{1+2\varepsilon}{d}\right)^{n},$$ 
for every $n\geq 0$.
Therefore 
$$\| (SG_{0})^{n+1} \| = \| SG_{0}(SG_{0})^n \| \leq \| S \|
M\left(\frac{1+2\varepsilon}{d}\right)^{n}.$$
We obtain $r(SG_{0}) \leq (1+2\varepsilon)/d$.
The fact that $G(\ld)$ is a generalized inverse of $T-\ld S$ implies 
$TG_{0}T = T$ and $G_{0}TG_{0} = G_{0}$. Hence
\begin{eqnarray*}
	\sup \{\frac{1}{r(SL)} : L \in \LL (Y,X), \;  
	TLT = T, LTL = L\} & \geq 
	&  \frac{1}{r(SG_{0})} \\
	 & \geq &  \frac{d}{1+2\varepsilon} .
\end{eqnarray*}
Since this holds for every $\varepsilon > 0$, 
we get that both suprema are not smaller than $d$.

For the second inequality, let $L\in 
\LL(Y,X)$ with $TLT = T$. Set 
$$F(\ld) = L\left( I - \ld SL\right)^{-1} = \sum_{k\geq 0}L(SL)^k\ld^k$$
which is defined and analytic for $\ld$ satisfying 
$$|\ld| < \alpha = \min (\|SL\|^{-1} ; \|LS\|^{-1} ).$$ 
By Lemma \ref{lem:2} we have
$$(T- \ld S) F(\ld)(T- \ld S) = T- \ld S$$
for each $\ld$ satisfying $|\ld| < \alpha$. Using Lemma 
\ref{lem:conorme}, this implies
$$\gm_{m}(T;S) \geq \frac{1}{\|L(SL)^{m-1}\|} 
\geq \frac{1}{\|L\|\; \|(SL)^{m-1}\|}.$$
Using \cite{BaLa}, we have
$$d = \lim_{m\to \infty} \gm_{m}(T;S)^{1/m} \geq 1/r(SL),$$
for every $L \in \LL (Y,X)$ satisfying $TLT = T$. This gives the 
desired inequality
$$
d \geq \sup \{\frac{1}{r(SL)} : L \in \LL (Y,X), \;  TLT = T\}.
$$
The proof is now complete. 
\myqed

\begin{rem}
	Using the Kato \cite{kato} 
	decomposition, the formula $d(T;S) = \lim_{n\to \infty} 
	\gm_{n}(T;S)^{1/n}$ was proved in \cite{BaLa} for $T$ Fredholm and 
	$S$ arbitrary, without the condition that $\dim N(T-\lambda S)$ is
constant in a neighborhood of $\lambda = 0$. 
	We do not know if Theorem \ref{thm:main} holds without this condition, equivalent
	to $k(T;S) = 0$.
\end{rem}

\noindent {\bf Acknowledgment.} We are grateful to the referee for 
valuable comments on this paper. 

\end{document}